\documentclass[12pt]{article}
\usepackage{latexsym}
\usepackage{amsfonts}
\usepackage{epsfig}
\newcommand{\binom}[2]{\left(\begin{array}{c}#1\\#2\end{array}\right)}

\def\sign{{\mathrm{sgn\,}}}
\def\inv{{\mathrm{inv\,}}}
\def\card{{\mathrm{card\,}}}
\def\ord{{\mathrm{ord\,}}}

\def\x{{\mathbf{x}}}
\def\y{{\mathbf{y}}}
\def\e{{\mathbf{e}}}

\def\Poly{{\mathrm{Poly}}}

\def\C{{\bf C}}    
\def\bC{{\bf\overline{C}}}
\def\R{{\bf R}}    

\def\N{{\bf N}}

\def\F{{\bf F}} 
\def\RP{{\bf RP}}

\def\FP{{\bf FP}}
\def\text{\hbox}   

\def\:{\colon}
\def\k{{\mathbf{k}}}

\def\q{{\mathbf{q}}}

\newtheorem{theorem}{Theorem}
\newtheorem{example}[theorem]{Example}
\newtheorem{corollary}[theorem]{Corollary}
\newtheorem{proposition}[theorem]{Proposition}
\newtheorem{lemma}[theorem]{Lemma}

\title{Degrees of real Wronski maps}
\author{A.\ Eremenko
 and A.\ Gabrielov\thanks{Both authors are
supported by NSF.}}
\date{March 21, 2002}
\begin{document}
\maketitle
\begin{abstract}
We study the map which sends vectors of polynomials
into their Wronski determinants. This defines a projection map of a 
Grassmann variety which we call a Wronski map.
Our main result is computation
of degrees of the {\em real} Wronski maps.
Connections with real algebraic geometry
and control theory are described.
\end{abstract}
\section{Introduction}

We study the map $W$ which sends vectors of polynomials
$(f_1,\ldots,f_p)$ to their
Wronski determinants:
\begin{equation}
\label{star}
W(f_1,\ldots,f_p)=\left|\begin{array}{ccc}f_1&\ldots&f_p\\
                                         f_1^\prime&\ldots&f_p^\prime\\
                                         \ldots&\ldots&\ldots\\
                                         f_1^{p-1}&\ldots&f_p^{p-1}\end{array}
\right|.
\end{equation}
Besides an intrinsic interest, this map is related to several
questions of algebraic geometry,
combinatorics and control theory as we describe below.

The following properties of the Wronski determinant
are well-known and easy to prove:
\newline
1. $W(f_1,\ldots,f_p)=0$ if and only if $f_1,\ldots,f_p$ are linearly
   dependent.
\newline
2. Multiplication of $(f_1,\ldots,f_p)$ by a constant matrix $A$
of size $p\times p$
results in multiplication of $W(f_1\ldots,f_p)$ by $\det A$.

These properties suggest that our map $W$ should be considered as
a map from a Grassmannian to a projective space. We recall the
relevant definitions.

Let $\F$ be one of the fields $\R$ (real numbers) or $\C$ (complex
numbers). For positive integers $m$ and $p$ we denote by
$G_\F=G_\F(m,m+p)$ the Grassmannian, that is the set of
all linear subspaces of dimension $m$ in $\F^{m+p}$. Such subspaces can be
described as row spaces of $m\times (m+p)$ matrices $K$ of maximal rank.
Two such matrices $K_1$ and $K_2$ define the same element of $G_\F$
if $K_1=UK_2$,
where $U\in GL(m,\F)$. It is easy to see that $G_\F(m,m+p)$
is an algebraic manifold
over $\F$ of dimension
$mp$. We have $G_\F(1,m+p)=\FP^{m+p-1}$,
the projective space over $\F$ of dimension $m+p-1$.

We may identify polynomials of degree at most $m+p-1$ in the domain
of the map $W$ in (\ref{star}) with vectors
in $\F^{m+p}$ using coefficients as coordinates, and similarly
polynomials in the range of $W$ with vectors in $\F^{mp+1}$.
Then, in view of the properties 1 and 2 of the Wronski determinant,
equation
(\ref{star}) will define a map $G_\F(p,m+p)\to\FP^{mp}$.
Alternatively, we can also identify polynomials of degree at most $m+p-1$
with linear forms on $\F^{m+p}$. Then $p$ linearly independent
forms define a subspace of dimension $m$ in $\F^{m+p}$, and we obtain a map
\begin{equation}
\label{wmap}
\phi:G_\F(m,m+p)\to\FP^{mp},
\end{equation}
which will be called a {\em Wronski map}.
To understand the nature of this map, we use the Pl\"ucker embedding
of the Grassmannian.

The Pl\"ucker coordinates of a point in $G_\F(m,m+p)$
represented by a matrix $K$
are the full size minors of $K$. This defines an embedding of $G_\F(m,m+p)$ to
$\FP^N,\;N=\binom{m+p}{m}-1$.
We usually identify $G_\F$ with its image under this embedding,
which is called a {\em Grassmann variety}.
It is a smooth algebraic variety in $\FP^N.$

Let $S\subset\FP^N$
be a projective subspace disjoint from $G_\F$, and $\dim_\F S=N-
\dim G_\F-1.$
We consider the {\em central projection} $\pi_S:\FP^N\backslash S\to\FP^{mp},$
and its restriction to $G_\F$,
\begin{equation}
\label{projection}
\phi_S=\pi_S\vert_{G_\F}: G_\F\to\FP^{mp}.
\end{equation}
Then $\phi_S$ is a finite regular map of
projective varieties. When $\F=\C$ this map has a degree,
which can be defined
in this case as the number of preimages of a generic point and is
independent of $S$.
This degree was computed by Schubert in 1886 (see \cite{K,GH,HP}
for modern treatment).
\vspace{.1in}

\noindent
{\bf Theorem A} {\em When $\F=\C$, the degree of $\phi_S$
is} 
\begin{equation}
\label{schubert}
d(m,p)=\frac{1!2!\ldots(p-1)!\,(mp)!}{m!(m+1)!\ldots(m+p-1)!}.
\end{equation}
\vspace{.1in}

Projective duality implies that $d(m,p)=d(p,m)$.
Here are some values of $d(m,p)$
$$
\begin{array}{ccccccccc}
&m=& 2    &3   &4     &5     &6      &7     &8\\
\\
p=2&& 2   &5   &14    &42    &132    &429     &1430\\
p=3&&     &42  &462   &6006  &87516  &1385670 &23371634\\
p=4&&     &    &24024 &1662804&140229804&\ldots&\ldots\\
p=5&&     &    &      &701149020&\ldots&\ldots&\ldots.
\end{array}
$$ 
In particular,
$$d(m,2)=\frac{1}{m+1}\binom{2m}{m},\quad\mbox{the $m$-th Catalan number.}$$
The numbers $d(m,p)$ have the following combinatorial interpretation:
they count the Standard Young Tableaux (SYT) of rectangular
shape $p\times m$. 

To see that the Wronski map is a projection
(\ref{projection})
we choose a center $S_0$ in the following way. Consider the
$p\times(m+p)$ matrix of polynomials
\begin{equation}
\label{normalcurve}
E(z)=\left(\begin{array}{c}F(z)\\F^\prime(z)\\ \ldots\\ F^{(p-1)}(z)\end{array}
\right),
\end{equation}
where $F(z)=\left(z^{m+p-1},\, z^{m+p-2}\,\ldots, z,\, 1\right)$.
For a fixed $z$, the row space of this matrix represents the osculating
$(p-1)$-subspace
to the rational normal curve $F:\FP^1\to\FP^{m+p-1}$ at the point $F(z)$.
The space $\Poly_\F^{mp}$ of
all non-zero polynomials $q\in\F[z]$ of degree at most $mp$,
up to proportionality,
will be identified with $\FP^{mp}$ (coefficients of polynomials serving as
homogeneous coordinates).

We claim that the Wrosnki map (\ref{wmap})
$\phi:G_\F\to\FP^{mp}$ can be defined by the formula
\begin{equation}
\label{map}
K\mapsto\phi(K)=
\det\left(\begin{array}{c}E(z)\\K\end{array}\right)\in\Poly_\F^{mp},
\end{equation}
where $K$ is a matrix of size $m\times(m+p)$ representing a point
in the Grassmannian $G_\F$.
First of all it is clear that (\ref{map}) indeed defines
a map $G_\F\to\FP^{mp}$: changing $K$ to $UK,
\; U\in GL(m,\F),$ will result
in multiplication of the polynomial $\phi(K)$ by $\det U$.
Furthermore, this map (\ref{map}), when expressed in terms of Pl\"ucker
coordinates, coincides with the restriction to $G_\F$ 
of
a projection of the form $\pi_S$ as in (\ref{projection}),
with some center which
we call $S_0$.
We do not need the explicit equations of $S_0$, but they can
be obtained by expanding the determinant in (\ref{map}) with respect to
the last $m$ rows, and collecting the terms with equal powers of $z$.

Now we verify that polynomial $\phi(K)$ in (\ref{map}) is a Wronski
determinant.
To see this, it is enough to consider the
``big cell'' $X$ of the Grassmannian $G_\F$, which is represented by the
matrices $K$ whose rightmost minor is different from zero.
We can normalize $K$ to make the rightmost $m\times m$ submatrix the unit
matrix.
If the remaining (leftmost) $p$ columns of $K$ are
$(k_{i,j}),\; 1\leq i\leq m,\; 1\leq j\leq p,$
then
$$\phi(K)=W(f_{1,K},\ldots,f_{p,K}),$$
where
\begin{equation}
\label{K12}
\begin{array}{lll}
f_{1,K}(z)&=&z^{m+p-1}-k_{1,1}z^{m-1}-\ldots-k_{m,1},\\
f_{2,K}(z)&=&z^{m+p-2}-k_{1,2}z^{m-1}-\ldots-k_{m,2},\\
\ldots&&\ldots\\
f_{p,K}(z)&=&z^{m}-k_{1,p}z^{m-1}-\ldots-k_{m,p}.
\end{array}
\end{equation}
Coefficients of these polynomials correspond to $p$ linear forms
that define the row space of the matrix $K=[(k_{ij}),I]$. This proves our
claim that (\ref{map}) coincides with the Wronski map.

For $p=2$, this
interpretation of the Wronski map as a projection
is due to L. Goldberg \cite{G}. Her notation for Catalan numbers
is different from our present notation.

In this paper we study the {\em real} map $\phi$, that is we set $\F=\R$.
One motivation of this study is the following conjecture due to 
B. and M. Shapiro: {\em If $w\in\Poly_\R^{mp}$ is a polynomial all of whose
roots are real, then the full preimage $\phi^{-1}(w)$ of this
polynomial consists of real points.} In \cite{EG1} we proved this conjecture
in the first non-trivial case $\min\{ m,p\}=2$. On the other hand,
when $m$ and $p$ are both even,
there are polynomials $w\in\Poly_\R^{mp}$ which do not have
real preimages under the Wronski map. So it was natural
to ask the question, whether for some $m$ and $p$ one can give a lower
estimate for the number of real preimages.
To asnwer this question, we compute in this paper the topological
degree of
the real Wronski maps.

Notice the following important property of $\phi$: it sends the big cell $X$
of the Grassmannian into
the big cell $Y$ of the projective space consisting of those polynomials
whose degree is exactly $mp$. Moreover, it sends the complement
$G_\F\backslash X$ into $\FP^{mp}\backslash Y$. When $\F=\R$ these cells
$X$ and $Y$ can be identified with $\R^{mp}$, in particular they are
orientable, and the restriction of
$\phi$ to $X$ is a smooth map
\begin{equation}
\label{restriction}
\phi:X\to Y,\quad \phi(\partial X)\subset\partial Y.
\end{equation}
To define the {\em degree} of such map (see, for example
\cite{Milnor}), we fix some orientations on $X$ and $Y$. Then choose a
regular value $y\in Y$ of $\phi$, which exists by Sard's theorem,
and define
\begin{equation}
\label{degree}
\deg \phi=\pm\sum_{x\in \phi^{-1}(y)}\sign\det \phi^\prime(x),
\end{equation}
using local coordinates in $X$ consistent with the chosen orientation of $X$,
and any local coordinate at $y$. 
The degree $\deg f$ changes sign if one changes one
of the orientations of $X$ or $Y$;
it is independent of the choice of
local coordinates within the class defined by the chosen
orientation of $X$, and of the regular
value $y$. 
In Section 3 we will discuss a more general definition of degree
which does not use special properties of the Wronski map and
applies to all equidimensional projections of real Grassmann
varieties.

%
%

To state the main result of this paper, we need a definition.
Consider the sequences $\sigma=(\sigma_j)$
of length $mp$ whose entries are elements of
the set $\{1,\ldots,p\}$, and each element occurs exactly $m$ times.
Suppose that the following additional condition is satisfied: for every
$n\in [1,mp]$ and every pair $i<k$ from $\{1,\ldots,p\}$, 
\begin{equation}
\label{ballot}
\#\{ j\in[1,n]:\sigma_j=i\}\geq\#\{j\in[1,n]:\sigma_j=k\}.
\end{equation}
Such sequences are called {\em ballot sequences}
or {\em lattice permutations}
\cite{MacMahon}.
For given $m$ and $p$, the set of all ballot
sequences is denoted by $\Sigma_{m,p}$.
There is a natural correspondence between $\Sigma_{m,p}$
and the set of the standard Young tableaux of rectangular shape $p\times m$
\cite[Proposition 7.10.3]{Stanley}: we fill the shape with integers from $1$
to $mp$ putting one integer in each cell;
if $\sigma_j=i$ we put the integer $j$
to the leftmost unoccupied place in the row $i$.
(As usual, the row number increases downwards).

Frobenius and
MacMahon independently
found that the cardinality of $\Sigma_{m,p}$ is $d(m,p)$,
the same number 
as in 
(\ref{schubert}), see for example, \cite[Sect. III, Ch. V, 103]{MacMahon}
or \cite[Proposition 7.21.6]{Stanley}. Of course, the coincidence of these 
numbers is not accidental \cite{Fulton,Stanley}.

Let $\sigma\in\Sigma_{m,p},\;\sigma=(\sigma_j)$. A pair $(\sigma_j,\sigma_k)$
is called an {\em inversion} if $j<k$ and $\sigma_j>\sigma_k$.
In terms of the SYT, and inversion occurs each time
when for a pair of integers the greater integer of the pair stands in
higher row than the smaller one. 
The total number of inversions in $\sigma$ is denoted by $\inv\sigma$.
Now we define
\begin{equation}
\label{seq}
I(m,p)=\left|\sum_{\sigma\in\Sigma_{m,p}}(-1)^{\inv\sigma}\right|.
\end{equation}
It is clear that $I(m,p)=I(p,m)$ because a pair of entries in
a SYT is an inversion if and only if the same pair in the transposed
SYT is not an inversion. This permits us to restrict to the case 
\begin{equation}
\label{wlog}
m\geq p\geq 2
\end{equation}
in the computation of the numbers $I(m,p)$. 
Recently, D. White \cite{White} proved that $I(m,p)=0$ iff $m+p$ is even.
For odd $m+p$ satisfying (\ref{wlog}),
he found that $I(m,p)$ coincides with the number of
shifted standard Young tableaux (SSYT) of shape
$$\left(\frac{m+p-1}{2},
\frac{m+p-3}{2}, \ldots, \frac{m-p+3}{2}, \frac{m-p+1}{2}
\right).$$
An explicit formula for the number of SSYT
(see, for example, \cite[Proposition 10.4]{Hoffman}) gives
$I(m,p)=$
$${1! 2! \cdots (p-1)! (m-1)! (m-2)! \cdots (m-p+1)! (mp/2)!
\over (m-p+2)! (m-p+4)!\cdots(m+p-2)!
\left(\frac{m-p+1}{2}\right)!\left(\frac{m-p+3}{2}\right)!\cdots
\left(\frac{m+p-1}{2}\right)!},$$
when $m+p$ is odd.
SSYT appear in Schur's theory of
projective representations of symmetric groups,
see, for example, \cite{Hoffman}.
Here are some values of $I(m,p)$: 
$$
\begin{array}{lcccccccccccc}
 m=&3&  4 & 5 &  6 &   7 &   8  &   9  &   10  &  11  &   12   \\ 
\\
p=2&1&  0 & 2 &  0 &   5 &   0  &  14  &  0    &  42  &    0  \\
p=3&0&  2 & 0 & 12 &0    & 110  &   0  &1274   &   0  & 17136\\
p=4& &  0 &12 & 0  &286  & 0    & 12376&  0    &759696&\ldots\\
p=5& &    & 0 & 286&0    &33592 &   0  &8320480&   0  &\ldots\\
\end{array}$$
When $p=2$ and $m$ is odd, $I(m,2)=d((m-1)/2,2)$, a Catalan
number. 
The main result of this paper is
\begin{theorem}\label{thm2}
The degree of the real Wronski map $(\ref{map})$ is
$\pm I(m,p)$. 
\end{theorem}
\begin{corollary}\label{cor1}
If $m+p$ is odd then the real Wronski map $(\ref{map})$ is surjective;
a generic point $y\in\RP^{mp}$ has at least
$I(m,p)$ real preimages.
\hfill$\Box$
\end{corollary}

It follows from Theorem A, that for all $y\in\Poly_\R^{mp}$,
\begin{equation}\label{card}
\card
\phi^{-1}(y)\cap G_\R(m,m+p)\leq d(m,p).
\end{equation}
This estimate
is best possible for every $m$ and $p$, see \cite{sottile}, or
the remark at the end of Section 2.

If both $p$ and $m$ are even, then Theorem \ref{thm2}
gives $\deg G_\R(m,m+p)=0$, and in fact in this case  
the preimage in Corollary \ref{cor1} may be empty,
as examples in \cite{EG4} show. 

The lower bound in Corollary \ref{cor1} with $p=2$ is 
best possible, as the following example given in \cite{EG3} shows:
\begin{example}\label{ex1} For $p=2$ and every odd $m$, there exist 
regular values $y\in\RP^{2m}$ 
such that 
the cardinality of $\phi^{-1}(y)$
is $I(m,2)=d((m-1)/2,2).$
\hfill$\Box$
\end{example}

To each $p$-vector $(f_1,\ldots,f_p)$ of linearly independent polynomials
one can associate
a rational
curve $f=(f_1:\ldots:f_p)$ in $\FP^{p-1}$, whose image is not contained
in any hyperplane. The following {\em equivalence relation}
on the set of rational curves
corresponds to the equivalence relation on the $p$-vectors of polynomials:
\begin{equation}
\label{equiv}
f\sim g\quad\mbox{if}\quad f=\ell\circ g,\quad
\mbox{where}\; \ell\;\mbox{is an automorphism of}\quad\FP^{p-1}.
\end{equation}

If $(f_1,\ldots,f_p)$ is a coprime $p$-vector of polynomials, then the roots of 
\newline
$W(f_1,\ldots,f_p)$ coincide with finite inflection points of the curve $f$.
Notice that $G_\R\subset G_\C$ can be represented by $p$-vectors
of real polynomials,
and to each such $p$-vector corresponds a real curve $f$.
When $p=2$, $f=f_2/f_1$ is a rational function. If the pair
$(f_1,f_2)$ is coprime, roots of $W(f_1,f_2)$
are the finite critical points of $f$. 
Thus our Theorem \ref{thm2} has the following
\begin{corollary}\label{cor7}
Let $X$ be a set of $mp$ points in general position
in $\bC$, symmetric with respect to $\R$.
Then the number $k$ of equivalence classes of real rational curves
in $\RP^{p-1}$ of degree
$m+p-1$ whose sets of inflection points coincide with $X$ satisfies $k\geq I(m,p)$.
In particular, for $p=2$, this number $k$
satisfies
\begin{eqnarray}
&0\leq k\leq d(m,2),\quad\mbox{if $m$ is even, and}\label{label1}\\   
&d((m-1)/2,2)\leq k\leq d(m,2),\quad\mbox{if $m$ is odd.}\label{label2}
\end{eqnarray}
\end{corollary}
Examples in \cite{EG3} show that for every $m$, the lower estimates in
(\ref{label1}) and (\ref{label2})
are best possible.
So when $m$ is odd, the Wronski map $\phi:G_\R(m,m+2)\to\RP^{2m}$
is surjective, while for even $m$ it is not.\hfill$\Box$

We prove Theorem \ref{thm2} in Section 2.
In Section 3 we discuss the definition of degree
for arbitrary projections of real Grassmann varieties and interpret
Theorem \ref{thm2} in terms of control theory. 

For the case $p=2$, the results of this paper were obtained in
\cite{EG3}, with a different method based on \cite{EG1}.
We thank
S. Fomin, Ch. Krattenthaler, F. Sottile and R. Stanley
for helpful suggestions. 

\section{Computation of degree}

In this section we prove Theorem \ref{thm2}.
We fix integers $m,p\geq 2$.
Consider vectors of integers $\k=(k_1,\ldots,k_p)$ satisfying
$$0\leq k_1<k_2<\ldots<k_p<m+p,$$ 
and vectors of real polynomials
$\q=(q_1,\ldots,q_p)$ of the form
\begin{equation}
\label{PQ}
\begin{array}{lll}
q_1(z)&=&z^m+a_{1,m-1}z^{m-1}+\ldots+a_{1,k_1}z^{k_1},\\
&&\\
q_2(z)&=&z^{m+1}+a_{2,m}z^m+\ldots+a_{2,k_2}z^{k_2},\\
\ldots&&\ldots\\
q_p(z)&=&z^{m+p-1}+a_{p,m+p-2}z^{m+p-2}+\ldots+a_{p,k_p}z^{k_p}.
\end{array}
\end{equation}
Suppose that all coefficients
$a_{ij},\; k_i\leq j\leq m+i-2,\; 1\leq i\leq p,$ are
positive, 
all roots of the Wronskian $W_\q=W(q_1,\ldots,q_p)$
belong to the semi-open interval
$(-1,0]\subset\R,$ and those roots on the open interval $(-1,0)$ are simple.
The set of all such polynomial vectors $\q$ will
be denoted by $b(\k)$.
The greatest common factor of $\{ q_1,\ldots,q_p\}$ is $z^{k_1}$.

It is easy to see that $b(\k)$ parametrizes a subset of the big cell
of the
Grassmannian $G_\R(m,m+p)$: the representation of a point of $G_\R(m,m+p)$
by a vector from $b(\k)$ is unique.
Setting $k_i=i-1$, $1\leq i\leq p$, we obtain an open subset
$b(0,1\ldots,p-1)\subset
G_\R(m,m+p)$. We define
$$k=k_1+(k_2-1)+\ldots+(k_p-p+1)\geq 0.$$
It is easy to see that $k$ is the multiplicity of the root of $W_\q$ at $0$
for $\q\in b(\k)$.
Using coefficients of $\q$ as coordinates,
we can identify $b(\k)$ with a subset
of $\R^{mp-k}$, and 
introduce an orientation by ordering these coefficients:
\begin{equation}
\label{order}
a_{1,k_1},\ldots,a_{1,m-1},a_{2,k_2},\ldots,a_{2,m},\ldots,a_{p,m+p-2}.
\end{equation}
In this sequence, coefficients
of $q_j$ precede coefficients of $q_k$ for $j<k$,
and coefficients of one polynomial $q_j$ are ordered according to their
second subscript. It is useful to place these coefficients into a
Young diagram $Y$ with $p$ rows,
such that coefficients of $q_i$ 
are in the $i$-th row, their second subscript decreasing left to right.
For $\q\in b(\k)$, we denote the negative roots of the Wronskian
$W=W_\q$ by
\begin{equation}
\label{zerosW}
-x_{mp-k}<-x_{mp-k-1}<\ldots<-x_1.
\end{equation}
In addition to these, there is a root of multiplicity $k$ at $0$.
We denote
by $\Delta_{\q}$ the Jacobi matrix of the map $b(\k)
\to\Poly_\R^{mp-k},\;
\q\mapsto W_\q$,
using coordinates (\ref{order}) in the domain and $\x=(x_1,x_2,\ldots,x_{mp-k})$
in the
range, where $-x_j$ are the negative roots of $W_\q$ as in (\ref{zerosW}).
So the $i$-th row of this matrix $\Delta_\q$ corresponds to $x_i$, and the 
$j$-th column to the $j$-th term of the sequence (\ref{order}). 
When $k=0$, so that $\k=(0,1,\ldots,p-1)$, and $b(\k)$ is an open subset
of $G_\R(m,m+p)$, we have $\Delta_{\q}=\phi'(\q)$, the
derivative of the Wronski map with respect to the chosen coordinates.

For example, $b(m-1,m+1,m+2\ldots,m+p-1)$ consists of vectors 
\begin{equation}
\label{u1}
q_1(z)=z^m+a_{1,m-1}z^{m-1},\quad q_2(z)=z^{m+1},\ldots, q_p(z)=z^{m+p-1},
\end{equation}
the Wronskian is
$$W(z)=(p-2)!z^{mp-1}\left((p-1)z+p!a_{1,m-1}\right),$$
and its only negative root is $-p!a_{1,m-1}/(p-1)$. So
\begin{equation}
\label{2geq0}
\det\Delta_{\q}=p!/(p-1)>0\quad\mbox{for}\quad \q\in b(m-1,m+1,\ldots,m+p-1).
\end{equation}
This example will be later used as a base of induction.

We denote by $E$ the set of all increasing homeomorphisms
$\epsilon:\R_{>0}\to\R_{>0},\; \epsilon(t)<t$ for $t>0$. Let $n\geq 1$
be an integer, and $\epsilon\in E.$
A {\em thorn} $T(n,\epsilon)$ in $\R^n$ is defined as
\begin{equation}
\{\x=(x_1,\ldots,x_n)\in\R^n_{>0}:\; x_j<\epsilon(x_{j+1}),\;
1\leq j\leq n-1,\; x_n<\epsilon(1)\}.\label{thorn}
\end{equation}
Notice that this definition depends on the ordering of coordinates in $\R^n$.
We always assume that this ordering corresponds to the increasing order
of subscripts.
\begin{lemma}
\label{l0}
Intersection of any finite set of thorns in $\R^n$ is a thorn in $\R^n$.
\end{lemma}
{\em Proof.} Take the minimum of their defining functions $\epsilon$.
\hfill$\Box$
\begin{lemma}
\label{l1} Let $T=T(n,\epsilon)$ be a thorn in $\R^n=\{(x_1,\ldots,x_n)\}$,
and $U$ its neighborhood
in $\R^{n+1}=\{(x_0,x_1,\ldots,x_n)\}.$
Then $U^+=U\cap\R^{n+1}_{>0}$ contains a thorn $T(n+1,\epsilon_1)$. 
\end{lemma}

{\em Proof.} There exists a continuous function $\delta:T\to\R_{>0}$,
such that $U^+$ contains the set $\{(x_0,\x):\x\in T, 0<x_0<\delta(\x)\}$.
Let $\delta_0(t)$ be the minimum of 
$\delta$ on the compact subset $\{\x\in \overline{T(n,\epsilon/2)}:x_1
\geq t
\}$
of $T$. Then there exists $\epsilon_0\in E$ with the property $\epsilon_0<
\delta_0$. 
If we define  
$\epsilon_1=\min\{\epsilon/2,\epsilon_0\}$, then $T(n+1,\epsilon_1)
\subset U^+$.
\hfill$\Box$
\begin{lemma}
\label{l2}
Let $T=T(n+1,\epsilon)$ be a thorn in $\R^{n+1}$, and $h:T\to\R^{n+1}_{>0},$
\newline
$(x_0,\x)\mapsto(y_0(x_0,\x),\y(x_0,\x))$,
a continuous map with the properties: for every $\x$ such that $(x_0,\x)\in T$
for some $x_0>0$,
the function $x_0\mapsto y_0(x_0,\x)$ is increasing,
and $\lim_{x_0\to 0}\y(x_0,\x)=\x$.
Then the image $h(T)$ contains a thorn.
\end{lemma}

{\em Proof}. We consider the region $D\in\R^{n+1}$ consisting of $T$,
its reflection $T'$ in the hyperplane $x_0=0$ and the interior
with respect to
this hyperplane 
of the common boundary of $T$ and $T'$.
The map $h$ extends to $T'$ by symmetry:
$h(-x_0,\x)=-h(x_0,\x), (x_0,\x)\in T$, and then to the whole $D$
by continuity. It is easy to see that the image of the extended map
contains a neighborhood $U$ of the 
intersection of $D$ with the hyperplane $x_0=0$. 
This intersection is a thorn $T_1$ in $\R^n=\{(x_0,\x)\in\R^{n+1}:
x_0=0\}$.
Applying Lemma
\ref{l1} to this thorn $T_1$, we conclude that $U^+$ contains a thorn.
\hfill$\Box$
\vspace{.1in}

Given an increasing homeomorphism $\epsilon\in E$,
we define 
$w(k,\epsilon)\subset \Poly_\R^{mp}$ as the set 
of 
all real monic polynomials of degree $mp$ with $mp-k$ negative roots as in
(\ref{zerosW}), these roots satisfying (\ref{thorn}) with $n=mp-k$,
and a root of
multiplicity $k$ at $0$. Thus $w(k,\epsilon)$ is parametrized by a
thorn $T(mp-k,\epsilon)$ in $\R^{mp-k}$.

Starting with $b(m-1,m+1,\ldots,m+p-1)$, 
we will generate subsets of 
$b(\k)$
by performing the following operations $F^i,\; 1\leq i\leq p$,
whenever they are defined.
Suppose that for some $i\in\{1,\ldots,p\}$ and some multiindex $\k$,
the following condition is satisfied:
\begin{equation}
\label{kcond}
i>1\quad\mbox{and}\quad k_i>k_{i-1}+1,\quad\mbox{or}\quad i=1\quad\mbox{and}
\quad k_1>0.
\end{equation}
Notice that for given $\k$, this condition is satisfied with {\em some} $i\in
\{1,\ldots,p\}$ iff $k>0$.
If (\ref{kcond}) holds,
we define a family of operators $F^i:b(\k)\to b(\k-\e_i)$, where
$\e_i$ is the $i$-th standard basis vector in $\R^p$, by
\begin{equation}
\label{oper1}
\q\mapsto F_a^i(\q)=(\q+az^{k_i-1}\e_i),
\end{equation}
where $a>0$ is a small parameter, whose range may depend on $\q$.
Thus an operation $F^i$ leaves all
polynomials in $\q$, except $q_i$, unchanged. 
The following Proposition shows, among other things,
that $F^i$ are well defined if the range of $a$ is appropriately
restricted. 
\begin{proposition}\label{Prop1}
Suppose that for some $\epsilon\in E$, and $\k$ and $i$
satisfying {\rm (\ref{kcond})},
a set $U\subset b(\k)$ is given,
such that the map
$\q\mapsto W_\q:U\to w(k,\epsilon)$ is surjective,
and 
\begin{equation}
\label{app}
\det\Delta_{\q}\neq 0\quad\mbox{for}\quad
\q\in U.
\end{equation} 
Then there exist $\epsilon^*\in E$
and a set $U^*\subset b(\k^*),$ where $\k^*=\k-\e_i$,
with the following properties.
Every $\q^*\in U^*$ has the form $F_a^i(\q)$ where $F_a^i$ is
defined in {\rm (\ref{oper1})}, $\q\in U$, and $a>0$;
\begin{equation}
\label{bi}
\mbox{the map}\quad \q^*\mapsto W_{\q^*}:U^*\to w(k-1,\epsilon^*)\quad
\mbox{is surjective,}
\end{equation}
and
$\det\Delta_{\q^*}\neq 0$ for $\q^*\in U^*$. Moreover,
\begin{equation}
\label{maineq}
\sign\det\Delta_{\q^*}=(-1)^{\chi(\k,i)}\,\sign\det\Delta_{\q},
\end{equation}
for every $\q^*\in U^*$ and every $\q\in U$, where
$\chi(\k,i)$ is the number of terms in the sequence {\rm (\ref{order})}
whose first subscript is less than $i$. In other words, $\chi(\k,i)$
is the total number of cells in the rows $1$ to $i-1$ in the Young diagram
$Y$ described after {\rm (\ref{order})}.
\end{proposition}
\vspace{.1in}

{\em Proof}. Let us fix $\q\in U$, and put $W=W_\q$. 
As $W\in w(k,\epsilon),$ we have $\ord W=k$, 
where $\ord$ denotes the multiplicity of a root at $0$.
Let $cz^k$ be the term of the smallest degree in $W(z)$.
Then $c>0$, because all roots of $W$ are non-positive.
In fact,
\begin{equation}
\label{c}
c=\prod_{j>l}(k_j-k_l)\prod_ja_{j,k_j}>0.
\end{equation}
We define $W^*=W_{\q^*},$ where $\q^*=F_a^i(\q)$.
Then $\ord W^*=k-1$ and the term of the smallest degree
in $W^*(z)$ is $c^*z^{k-1}$, where
\begin{equation}
\label{cstar}
c^*=a\prod_{j>l}(k^*_j-k^*_l)\prod_{j\neq i}a_{j,k_j}>0.
\end{equation}
We conclude that when $a$ is small enough (depending on $\q$), the
Wronskian $W^*$ has one simple root in a neighborhood of each negative root
of $W$, and in addition, one simple negative root close to zero, and
a root of multiplicity $k-1$ at $0$.
To make this more precise, we denote the negative roots of $W$ and $W^*$ by
\begin{equation}
\label{igrek}
-x_n<\ldots<-x_1\quad\mbox{and}\quad
-y_{n}<\ldots<-y_{1}<-y_0,
\end{equation}
where $n=2m-k$, and $y_j=y_j(a)$.
We have 
\begin{equation}
\label{trah}
y_j(0)=x_j,\quad\mbox{for}\quad 1\leq j\leq n,\quad\mbox{and}\quad y_0(0)=0.
\end{equation}
Furthermore, if $a$ is small enough (depending on $\q$)
\begin{equation}
\label{trah2}
a\mapsto y_0(a)\quad\mbox{is increasing and continuous}.
\end{equation}

The set $w(k,\epsilon)$ is parametrized by a thorn
$T=T(n,\epsilon)$,
where $\x=(x_1,\ldots,x_n)$, and $n=mp-k$.
There exists a continuous function $\delta_0:T\to\R_{>0}$, such that
\begin{equation}
\label{11}
\q^*\in b(\k^*),\quad\mbox{for}\quad a\in (0,\delta_0(\x)),\quad
\x\in T.
\end{equation}
Now we are going to compare $\det\Delta_{\q}$ with $\det\Delta_{\q^*}$.
For this purpose we investigate the asymptotic behavior of the `new root'
$y_0(a)$ of $W^*$, as $a\to 0$.
Comparison of the terms of the lowest degrees in $W(z)$ and $W^*(z)$,
(\ref{c}) and (\ref{cstar}) show
that
\begin{equation}
\label{behaves1}
-y_0(a)=-c^*/c+o(a)=-c_\k a/a_{i,k_i}+o(a),\quad a\to 0,
\end{equation}
where $c_\k>0$ depends only on the multiindex $\k$.

The Jacobi matrix $\Delta^*=\Delta_{\q^*}$ is obtained
from the Jacobi matrix $\Delta=\Delta_{\q}$ by adding the top
row,
corresponding to $y_0$, and a column, corresponding to $a_{i,k_i-1}=a$.
The position of the added column is 
$$1+\chi(\k,i),$$
where $\chi(\k,i)$ is the number of terms of the sequence (\ref{order})
whose first subscript is less than $i$.

According to (\ref{behaves1}),
the intersection
of the added row with the added column contains the only essential element
of this row: 
$$\partial y_0/\partial a=c_\k/a_{i,k_i}+o(1),\quad a\to 0.$$
The rest of the elements of the
first row of $\Delta^*$ are $o(1)$ as $a\to 0$.
Expanding $\Delta^*$ with respect to its first row, we obtain
$$\det\Delta^*=(-1)^{\chi(\k,i)}\, (c_\k/a_{i,k_i})\det\Delta+o(1),\quad a\to 0.$$
Now it follows from our assumption (\ref{app}) that for sufficiently
small $a$,  $\Delta^*\neq 0$. Moreover, (\ref{maineq}) holds, if $a$ is sufficiently
small. 
More precisely,
for every $\q\in U$ there exists
$\delta_1(\q)>0$ such that for $0<a<\delta_1(\q)$ we have
$\det\Delta^*\neq 0$, and (\ref{maineq}). 
Taking $\delta=\min\{\delta_0,\delta_1\}$, where $\delta_0$ was defined in
(\ref{11}), we obtain the set
\begin{equation}
\label{neighbor}
U^*=\{\q^*=F_a(\q_\x):\x\in T, a\in (0,\delta(\x))\}\subset b(\k^*),
\end{equation}
and this set $U^*$ satisfies (\ref{maineq}). Here $\q_\x\in U$ is
some preimage of $W_\x\in w(k,\epsilon)\cong T$ under the map
$\q\to W_\q$. Such preimage exists
by assumption of Proposition \ref{Prop1} that the map
$\q\mapsto W_\q,$  $U\to w(k,\epsilon)$ is surjective. It remains to achieve
(\ref{bi}) by modifying the thorn $T$. This we do in two steps.
First we
apply Lemma \ref{l1} to the half-neighborhood (\ref{neighbor}) of
$T$, with $x_0=a$, to obtain a thorn $T_1(n+1,\epsilon_1)$ in $\R^{n+1}$.
Then we apply Lemma \ref{l2} to the map $h:T_1\to\R^{n+1}_{>0}$,
defined by $y_j=y_j(x_0,\x)$, where $y_j$ are as in (\ref{igrek}), and $x_0=a$.  
This map $h$ satisfies all conditions of Lemma \ref{l2} in view of
(\ref{trah}) and (\ref{trah2}). This proves (\ref{bi}).
\hfill$\Box$
\vspace{.1in}

{\em Conclusion of the proof of Theorem \ref{thm2}}.
We begin with a brief outline of our argument.
For each $\k$ and $i$ satisfying (\ref{kcond}), equation 
(\ref{oper1}) defines an operator depending on parameter $a$:
$F_a^i:b(\k)\to b(\k-\e_i).$ 
Starting from a subset of
$b(m-1,m+1,\ldots,m+p-1)$, we will consecutively apply operators $F^i$ 
in all possible sequences allowed by (\ref{kcond}).
 In the end we will obtain a set of polynomial $p$-vectors
in $b(0,1,\ldots,p-1)$, which will contain the full preimage of a point under
the Wronski map. Equations (\ref{maineq}) will permit
to control the sign of the Jacobian determinant of the Wronski map at all
points of this preimage.

Now we give the details.
Consider the set of all finite (non-empty)
sequences $\sigma=(\sigma_j)$,
where $\sigma_j\in\{1,\ldots,p\},\; j\in \N,$ satisfying 
(\ref{ballot}). For every such sequence we define
$$\k(\sigma)=(k_1,\ldots,k_p),\quad\mbox{where}\quad k_i=m+i-1-\#\{ j:
\sigma_j=i\},$$
and 
$$k(\sigma)=\sum_{i=1}^p k_i(\sigma)-i+1.$$
Let $\Sigma=\Sigma(m,p)$ be the set of all sequences $\sigma$
satisfying (\ref{ballot}) and
$k(\sigma)\geq 0.$
Notice that for $\k=\k(\sigma)$, condition (\ref{kcond}) holds with some
$i\in\{1,\ldots,p\}$ if and only if $k(\sigma)>0$.

To each sequence $\sigma\in\Sigma$ we put into correspondence an open
set $U_\sigma\subset b(\k(\sigma))$ in the following way.
For $\sigma=(1)$,
we set 
$$U_{(1)}=\{\q\in b(m-1,m+1,\ldots,p-1): W_\q\in w(1,\epsilon_0)\},$$
where $\epsilon_0(x)=x$.
Then $U_{(1)}$ consists of the polynomial vectors of the form (\ref{u1})
with 
and $a_{1,m-1}\in (0,(p-1)/p!)$.

Applying operations $F^i$ to $U_{(1)}$ means that we use
Proposition \ref{Prop1} with
$U=U_{(1)}$, and $\k=(m-1,m+1,\ldots,p-1)$. We obtain from this Proposition the
sets $U^*$, which we call $U_{(1,i)}$.
In fact, Proposition \ref{Prop1} can we applied in this situation only with
$i=1$ or $i=2$.
Then we apply operations $F^j,\; j\in \{1,\ldots,p\}$ to
$U_{(1,i)}$, whenever permitted by (\ref{kcond}) and so on.

In general, suppose $U_\sigma$ is already constructed.
If $\k(\sigma)$ and $i$ satisfy (\ref{kcond}), we apply
operation $F^i$ to $U_\sigma$.
This means that we use Proposition \ref{Prop1} with
$U=U_\sigma$,
$\k=\k(\sigma)$ and this $i$.
The resulting $U^*$ is called $U_{(\sigma,i)}\subset b(\k(\sigma)-\e_i)$.

Every sequence $\sigma\in\Sigma$ 
encodes an admissible sequence of applications of operations $F^i$.
If $\sigma=(\sigma_j)$, then $\sigma_j=i$ indicates that
$F^i$ was applied on the $j$-th step.
Conditions (\ref{ballot}) and $k(\sigma)>0$ imply (\ref{kcond}) with some $i$,
so that an operation
$F^i$ is 
applicable. Every operation decreases $k(\sigma)$ by $1$, so the procedure
stops when $k(\sigma)=0$. 

Proposition \ref{Prop1} implies that for each $\sigma\in\Sigma$ with
$k(\sigma)\geq 0$,
there exists $\epsilon_\sigma\in E$ such that
\begin{equation}
\label{unramified}
W:U_\sigma\to w(k(\sigma),\epsilon_\sigma)
\end{equation}
is surjective and unramified. 

Observe that we can always replace $\epsilon^*$ in Proposition \ref{Prop1}
by a smaller function from the set $E$.
We use this observation
to arrange
that the coefficient, added to polynomials in $\q$ on each step,
is strictly smaller than all coefficients added on the previous steps.
This implies that for each $\q\in U_\sigma$, all coefficients
are strictly ordered, and the sequence $\sigma$ can be recovered from
this order. More precisely, let $k=k(\sigma)$, 
and $c_1>c_2>\ldots>c_{2m-k}>0$ be the ordering of the
sequence of coefficients of $q_1\ldots,q_p$.
Then $\sigma_j=i$ if $c_j=a_{i,l}$ with some $l$. In other words,
enumerating the cells of the
Young diagram $Y$ defined after (\ref{order}) in the order
of decrease of their entries gives a standard Young tableau.
The sequence $\sigma$ can be recovered from this tableau
in a unique way. 

We recall that the number
of inversions $\inv\sigma$ was defined in the introduction,
just before the equation (\ref{seq}).
We claim that for every $\sigma\in\Sigma$,
\begin{equation}\label{main2}
\sign\det \Delta_\q=\mu(\sigma)(-1)^{\inv\sigma}
\quad\mbox{if}\quad\q\in U_\sigma, 
\end{equation}
where $\mu(\sigma)=\pm1$ depends only on the length of $\sigma$.
Indeed, by (\ref{maineq}), on each step the sign of $\det\Delta$
is multiplied by $(-1)^{\chi}$, where $\chi=\chi(\k(\sigma),i)$
is the number of terms of $\sigma$ which are less than $i$.
This proves (\ref{main2}).

Now we consider the subset 
$$\Sigma_{m,p}=\{\sigma\in\Sigma_m:k(\sigma)=0\}.$$
It consists of ballot sequences, as defined in the introduction.
The set $\Sigma_{m,p}$ corresponds to rectangular standard Young tableaux
of the shape $p\times m$.
The number of such tableaux is
$d(m,p)$ (see, for example, \cite[Proposition 7.21.6]{Stanley}).
Sequences $\sigma\in\Sigma_{m,p}$ generate $d(m,p)$ open sets $U_\sigma
\subset b(0,1,\ldots,p-1)$
with the property that the maps (\ref{unramified}) are surjective
and unramified. Using Lemma \ref{l0}, we restrict these maps so that
they have a common range $w(0,\epsilon)$ with some $\epsilon\in E$.

As all maps (\ref{unramified}) are surjective,
every point from this common range has at least one preimage under
the Wronski map in each $U_\sigma, \sigma\in\Sigma_{m,p}$. All these
$d(m,p)$ preimages are different as elements of $b(0,1,\ldots,p-1)$, because
the sequence $\sigma$ can be recovered from the ordered
sequence of coefficients of $\q\in U_\sigma$. Furthermore,
all these $d(m,p)$ polynomial vectors represent different points in
the Grassmannian $G_\R(m,m+p)$, because to each point in $b(0,1,\ldots,p-1)$
corresponds
only one point of $G_\R(m,m+p)$.
Thus we found $d(m,p)$ different preimages of a point under the Wronski map.
On the other hand, the complex Wronski map has degree $d(m,p)$
by Theorem A, so we found all preimages of the real or complex Wronski map.
Equation (\ref{main2}) gives the signs of Jacobian determinants at these
points, so the degree of the Wronski map is given by (\ref{seq})
\hfill$\Box$
\vspace{.1in}

{\em Remark}. In the process of this proof, we constructed a point
in $\RP^{mp}$ which has $d(m,p)$ distinct real preimages under the Wronski map.
This proves the fact earlier established by Sottile \cite{sottile},
that the upper estimate $d(m,p)$ given by (\ref{card}),
is best possible for every $m$ and $p$.  

\section{Additional comments} 

\noindent
1. Let us show how to define topological degree (an unsigned integer)
for arbitrary projections
of real Grassmann varieties as in (\ref{projection}) with $\F=\R$.

Let $f:X\to Y$ be a smooth map of compact, connected
real manifolds of equal
dimensions. If $X$ is orientable, the degree $\deg f$ can be defined
by formula (\ref{degree}). If $X$ is orientable
but $Y$ is not then $\deg f=0$.

Now we suppose that both $X$ and $Y$ are non-orientable and consider
canonical orientable $2$-to-$1$ coverings
$\widetilde{X}\to X$ and $\widetilde{Y}\to Y$, which are
called the {\em spaces of orientations} of $X$ and $Y$ \cite[10.2]{Bi}.
The set $\widetilde{X}$ consists of pairs $(x,O)$ where $x\in X$ and
$O$ is one of the two orientations of the tangent space $T_x$.
There is a unique structure of smooth manifold on $X$ which makes
the map $\widetilde{X}\to X,\; (x,O)\mapsto x$ a covering, and $O$
depends continuously on $x$. The group of the covering 
$\widetilde{X}\to X$ is $\{\pm1\}$. Notice that the spaces
$\widetilde{X}$ and $\widetilde{Y}$ have canonical orientations.
A map $f:X\to Y$ is called {\em orientable} if there
exists a lifting $\widetilde{f}:\widetilde{X}\to\widetilde{Y}$, which commutes
with the action of $\{\pm1\}$, \cite[(10.2.5)]{Bi}.
A different but equivalent definition of an orientable map
is given in \cite[\S 5]{derham}.

For an orientable map, we
define $\deg f:=\pm\deg\widetilde{f}$.
Under our assumption that $X$ and $Y$ are connected,
this degree is defined
up to sign, which depends on the choice of the lifting.
Though $Y$ is connected,
$\widetilde{Y}$ may consist of one or two components, but the degree
is independent of the choice of the regular value $y\in\widetilde{Y}$.

Suppose that there exists a regular value $y\in Y$ and an affine chart
$U\subset X$,
so that $f^{-1}(y)\subset U$. Then we can compute the sum (\ref{degree})
using coordinates in $U$. The orientability of $f$ ensures that
this sum is independent of the choice of the chart $U$ and coincides
with $\deg f$.

To apply this construction to projections of Grassmann varieties
$$\phi_S:G_\R(m,m+p)\to\RP^{mp},$$
we recall that $\RP^{mp}$ is orientable iff $mp$ is odd,
and $G_\R(m,m+p)$ is orientable iff $m+p$ is even
(see, for example, \cite[Ch. 3 \S 2]{FR}). So in the case that $m+p$ is
even,
the degree of $\phi_S$ is defined
in the usual sense, as in (\ref{degree}).

To deal with the case when $m+p$ is odd, we first identify
the space of orientations of a Grassmannian $G_\R(m,n)$, with odd $n$.
Consider the
``upper Grassmannian''
${G}_\R^+(m,n)$,
which consists of all oriented $m$-subspaces in $\R^n$. It can be also
described as
the set of all $m\times n$ matrices $K$
of maximal rank, modulo the following equivalence relation:
$K'\sim K$ if $K'=UK$, where $\det U>0$. We have the natural $2$-to-$1$
covering ${G}_\R^+(m,n)\to G_\R(m,n)$, which assigns to the class
of $K$ in ${G}_\R^+(m,n)$ the class of the same $K$ in $ G_\R(m,n)$.
We also have
${G}_\R^+(1,n)=({\RP}^{n-1})^+$,
a sphere of dimension $n-1$.

Every upper Grassmannian is orientable and has canonical orientation.
To see this, we consider
the tangent space $T_x=T_x({G}_\R^+(m,n))$ which is the product of
$m$ copies of a subspace $y\cong\R^{n-m}$, complementary to $x$.
Orientation of $x$ induces a unique orientation of $y$, such that
$x\oplus y\cong\R^{n}$ has the standard orientation. This defines an orientation
on each tangent space $T_x$ which varies continuously with $x$.
So we have a  canonical orientation of ${G}_\R^+(m,n)$.

We claim that for odd $n$, the coverings
\begin{equation}
\label{cove}
G^+_\R(m,n)\to G_\R(m,n)\;\; \mbox{and}\;\; \widetilde{G}_\R(m,n)\to
G_\R(m,n)\quad\mbox{are isomorphic}.
\end{equation}
Indeed, for $x\in G_\R(m,n)$,
orientation of $x\subset\R^n$ defines an orientation
of $T_x(G_\R(m,n))$, as explained above. One can easily show that (in the
case of odd $n$) changing the
orientation of $x$ changes the orientation of $T_x(G_\R(m,n))$. This proves
(\ref{cove}).
\vspace{.1in}

We recall that a projection map $\pi_S:\RP^N\backslash S\to\RP^k$
can be described in homogeneous coordinates as
\begin{equation}
\label{matrix}
y=Ax,
\end{equation}
where
$A$ is a $(k+1)\times (N+1)$ matrix of maximal rank, and $x,y$ are
column vectors of homogeneous coordinates. The null space of $A$
represents the center of projection $S=S(A)$ (where the map is undefined).
Two matrices define the same projection if they are proportional.
A change of homogeneous coordinates in $\RP^N$ or in the target space $\RP^k$
results in multiplication of $A$
by a non-degenerate matrix from the right or left, respectively.
\begin{proposition}\label{propA} Let $n$ be an odd integer,
$G_\R(m,n)\subset\RP^N$ a Grassmann
variety, and $\phi:G_\R(m,n)\to\RP^{m(n-m)}$ a central projection.
Then $\phi$ is orientable.
\end{proposition}

{\em Proof}.
The Pl\"ucker embedding ${\mathrm{Pl}}:G_\R(m,n)\to\RP^N$ lifts to
$G^+_\R(m,n)\to(\RP^N)^+$, which is defined by the same rule as ${\mathrm{Pl}}$.
Using (\ref{cove}), we identify $\widetilde{G}_\R(m,n)$ with $G^+_\R(m,n)$,
and obtain the lifting
\begin{equation}\label{lift0}
\widetilde{{\mathrm{Pl}}}:\widetilde{G}_\R(m,n)\to(\RP^N)^+
\end{equation}
of ${\mathrm{Pl}}$.
Suppose now that a projection $\pi_S$ is defined by (\ref{matrix}) where
$A$ is an
\newline
$(m(n-m)+1)\times (N+1)$ matrix.
Then the same equation (\ref{matrix}) defines the lifting
\begin{equation}
\label{lift}
{\pi}^+:(\RP^N)^+\backslash{S}^+\to(\RP^{m(n-m)})^+\cong
\widetilde{\RP}^{m(n-m)},
\end{equation}
where $S^+$ is the preimage of $S$ under the covering
$(\RP^N)^+\to\RP^N,$ and the last isomorphism holds because
$m(n-m)$ is even.
Composition of the maps (\ref{lift0}) and (\ref{lift}) is the desired lifting
of $\phi_S$, which is evidently compatible with the action of $\{\pm1\}$.
The existence of such a lifting proves that $\phi_S$ is orientable.
\hfill$\Box$
\vspace{.1in}

Thus equidimensional
projections of real Grassmann varieties always have
well-defined degrees.
It is clear that when the center of projection
varies continuously, the degree does not change until the center $S$
intersects the Grassmann variety. These exceptional centers form
a subvariety $Z$ of codimension $1$ in the Grassmannian
$G_\R(N+1,N-mp)$ of all centers. So the degree is constant on
every component of $G_\R(N+1,N-mp)\backslash Z$, in particular,
all projections $\phi_S$ whose centers $S$ belong to the same component
of  $G_\R(N+1,N-mp)\backslash Z$
as the center $S_0$ of the Wronski map have the same degree $\pm I(m,p)$.
\vspace{.1in}

\noindent
2. Now we restate our results in terms of control theory by static output
feedback.
Suppose that a triple of real matrices
$\Sigma=(A,B,C)$ of sizes $n\times n$,
$n\times m$ and
$p\times n$ is given. This triple $\Sigma$ defines a
{\em linear system}
\begin{equation}
\label{system}
\begin{array}{lll}\dot{x}&=&Ax+Bu,\\ y&=&Cx.\end{array}
\end{equation}
Here $x,u$ and $y$ are functions of time (a real variable) taking their values
in $\R^n$, $\R^m$ and $\R^p$, respectively.
The values of these functions at a point $t\in\R$ are interpreted as
the state, input and output of our system at the moment $t$.

Behavior of the system (\ref{system}) is completely
determined by its 
{\em transfer function} $z\mapsto C(zI-A)^{-1}B$, which is a function of a complex
variable $z$ with values in the set of $p\times m$ matrices.
One wishes to control 
a given system (\ref{system}) by arranging a feedback,
which means sending the output to the input
via an $m\times p$ matrix $K$, called a {\em gain} matrix:
\begin{equation}
\label{gain}
u=Ky.
\end{equation}
Elimination of $u$ and  $y$ from (\ref{system}), (\ref{gain})
gives the {\em closed loop system}
$$\dot{x}=(A-BKC)x,$$
whose transfer function has poles at the zeros of the polynomial
\begin{equation}
\label{ppmap}
\psi_K(z)=\det(zI-A-BKC).
\end{equation}
The map $K\mapsto\psi_K\in\Poly_\R^{n}$ is called the {\em pole placement map},
and the problem of pole assignment is: {\em given a system $\Sigma$, and a 
set $\{ z_1,\ldots,z_n\}$,
symmetric with respect to $\R$, to find a real gain matrix $K$, such that
the zeros of $\psi_K$ are $\{ z_1,\ldots,z_n\}$.} Thus for a fixed system
$\Sigma$, arbitrary pole assignment is possible iff the pole
placement map is surjective. 

When $n>mp$, X. Wang \cite{wang} proved that for generic $\Sigma$,
the pole placement map is surjective.
Here we consider the case
$n=mp$.
We also assume that $n$ is the smallest possible size of a matrix $A$
in the representation $C(zI-A)^{-1}B$ of the transfer function.
Systems with this property are called ``controllable and observable'',
and and they form an open dense subset of the set
of all systems with fixed $(m,n,p)$.
To understand the structure of the pole placement map, we use
a coprime factorization of the open loop transfer function of a generic system
$\Sigma$
(see, for example, \cite[Assertion 22.6]{Delchamps}):
\begin{equation}
\label{dop}
\end{equation}
$$C(zI-A)^{-1}B=D(z)^{-1}N(z),\quad \det D(z)=\det(zI-A),$$
where $D$ and $N$ are polynomial matrix-functions of sizes 
$p\times p$ and $p\times m$, respectively. The polynomial matrix
$[D(z),N(z)]$ has the following properties: its full size minors
have no common zeros, and exactly one of these minors, $\det D(z)$,
has degree $n$ while all other minors have strictly smaller degree.
Every $p\times (m+p)$ polynomial matrix with these properties is related
to a linear system of the form (\ref{system}) via equations (\ref{dop}).

Using the factorization (\ref{dop}) and
the identity 
$\det(I-PQ)=\det(I-QP)$, which is true for all rectangular matrices of
appropriate dimensions,
we write
\begin{eqnarray*}
\psi_K(z)&=&\det(zI-A-BKC)=\det(zI-A)\det(I-(zI-A)^{-1}BKC)\\
&=&
\det(zI-A)\det(I-C(zI-A)^{-1}BK)\\
&=&\det D(z)\det(I-D(z)^{-1}N(z)K)
=\det(D(z)-N(z)K).
\end{eqnarray*}
This can be rewritten as
\begin{equation}
\label{12}
\psi_K(z)=\left|\begin{array}{cc}D(z)&N(z)\\
                                 K&I\end{array}\right|\in\Poly_\R^{mp}.
\end{equation}
In the last determinant, the first $p$ rows depend only on the 
given system, and the last $m$ rows on the gain matrix. 
Permitting arbitrary
$m\times(m+p)$ matrices $\hat{K}$ of maximal rank
as the last $m$ rows of the determinant in
(\ref{12}) we extend the pole placement map to
\begin{equation}
\label{extension}
\phi_\Sigma:G_\R(m,m+p)\to \RP^{mp},\quad \phi_{\Sigma}(\hat{K})=[\psi_K],
\end{equation}
where $[.]$ means the class of proportionality of a polynomial, which is
identified with a point in $\RP^{mp}$,
using the coefficients of a polynomial as homogeneous coordinates.
The map (\ref{extension}) is defined if for every matrix $\hat{K}$ of
rank $m$ in the last $m$ rows of (\ref{12}) the determinant in (\ref{12})
does not vanish identically. Systems $\Sigma$ with this property are called
{\em non-degenerate}, and they form an open dense subset in the set
of all systems with given $(m,p)$ and $n=mp$. 
Applying Laplace's expansion along the first $p$ rows
to the determinant in (\ref{12}), we conclude that the map $\phi_\Sigma$,
when expressed in 
Pl\"ucker coordinates, is nothing but a projection of
the Grassmann variety $G_\R(m,m+p)$ into $\RP^{mp}$
from some center depending on
$\Sigma$.
This interpretation of the pole placement map as a projection comes from
\cite{wang}. Now we notice that all projections arising from linear systems
as in (\ref{12}) have the property that they send the big cell $X$ of 
$G_\R(m,m+p)$ represented by matrices $\hat K$ of the form $[K,I]$
into the big cell $Y$ of $\Poly_\R^{mp}$ consisting of polynomials
of exact degree $mp$. Furthermore,
$G_\R(m,m+p)\backslash X$ corresponds to $\Poly_\R^{mp}\backslash Y$ inder such
projections. 
Our arguments in the first part of this section imply
that the real
pole placement maps of a non-degenerate system has a well-defined
degree\footnote{Projections arising from the linear systems 
map a fixed big cell $X$ of the Grassmannian to a fixed big cell
$Y$ of the projective space, and also send $\partial X$ to $\partial Y$.
This permits to define the degrees of these projections as we did it
for the Wronski map in the Introduction}.
As the center of projection $S=S(\Sigma)$
varies continuously, this degree remains
constant as long as $S$ does not intersect
the Grassmann variety. Degenerate systems are precisely those for which $S$
intersects the Grassmann variety.

Comparing (\ref{map}) with (\ref{extension}) we conclude
that the Wronski map is a pole placement map for some special
linear system.
So our Corollary \ref{cor1} implies
\begin{corollary}
\label{cor17}
For every $m$ and $p$ such that $m+p$ is odd
there is an open set $U$ of linear systems with $m$ 
inputs, $p$ outputs and state of dimension
$mp$, such that  for systems in $U$
the real pole placement map is surjective.
Furthermore, the pole placement problem for systems in $U$ has at least
$I(m,p)$ real solutions for any generic set of $mp$ poles
symmetric with respect to the real line.\hfill$\Box$
\end{corollary}

\vspace{.1in}

Purdue University, West Lafayette, Indiana 47907

eremenko@math.purdue.edu

agabriel@math.purdue.edu

\begin{thebibliography}{1}
\bibitem{Bi} N. Bourbaki, Vari\'et\'es diff\'erentielles et analytiques.
Fasc. de r\'esultats, Paragraphes 8 \`a 15, Hermann, Paris, 1971.
\bibitem{Delchamps}  D. Delchamps, State space and input-output linear systems,
Springer, NY, 1988.
\bibitem{EG1} A. Eremenko and A. Gabrielov, Rational functions with real
critical points and the B. and M. Shapiro conjecture in real enumerative
geometry, Ann. Math., 155 (2002) 109--129..
\bibitem{EG3} A. Eremenko and A. Gabrielov, The Wronski map and Grassmannians
of real codimension $2$-subspaces, Computational Methods and Function
Theory, to appear in vol. 1.
\bibitem{EG4} A. Eremenko and A. Gabrielov, Pole placement by static output
feedback for generic linear systems, to appear in SIAM Journal on
Control.
\bibitem{FR} D. Fuks and V. Rokhlin,
Beginner's course in topology. Geometric chapters.
(Transl. from Russian)
Springer-Verlag, Berlin-New York, 1984.
\bibitem{Fulton} W. Fulton, Young Tableaux, Cambridge UP, 1997.
\bibitem{G} L. Goldberg,
Catalan numbers and ramified coverings of the sphere,
Adv. Math.,  85 (1991) 129-144.
\bibitem{GH} Ph. Griffiths and J. Harris, Principles of algebraic geometry,
Willey, NY, 1978.
\bibitem{HP} W. Hodge and D. Pedoe, Methods of algebraic geometry, v. 2,
Cambridge UP,
1953.
\bibitem{Hoffman} P. Hoffman and J. Humphreys, Projective representations
of the symmetric groups, Clarendon Press, Oxford, 1992.
\bibitem{K} S. Kleiman, Problem 15. Rigorous foundation of Schubert's
enumerative calculus, in: F. Browder, ed.,
Mathematical development arising from Hilbert problems,
Proc. Symp. Pure Math., v. 28, AMS, Providence, 1976.
\bibitem{MacMahon} P. MacMahon, Combinatory Analysis, 2 vols.,
Cambridge UP, 1915-16 (reprinted by Chelsea, NY, in one vol. in 1960).
\bibitem{Milnor} J. Milnor, Topology from the differentiable viewpoint,
UP Virginia, 1965.
\bibitem{derham} G. de Rham, Differentiable manifolds, (English transl.)
Springer, NY 1984.
\bibitem{sottile} F. Sottile, Special Schubert calculus is real, Electronic Res.
Announcements, AMS, 5 (1999) 35-39.
\bibitem{Stanley} R. Stanley,  Enumerative Combinatorics, Vol. 2,
Cambridge Univ. Press, Cambridge, 1999.
\bibitem{wang} X. Wang,
Grassmannian, central projection,
and output feedback pole assignment of linear systems.
IEEE Trans. Automat. Control 41 (1996), no. 6, 786--794. 
\bibitem{White} D. White, Sign-balanced posets, J. Combin. Theory,
95 (2001), 1-38.
\end{thebibliography}
\end{document}